\documentclass[twocolumn]{autart}    

\usepackage{natbib}
\usepackage{amsfonts,amsmath,latexsym}
\usepackage{graphicx,xcolor,mathrsfs,mathtools}

\usepackage{balance}     
\usepackage{hyperref}

\newtheorem{definition}{Definition}
\newtheorem{proposition}{Proposition}
\newtheorem{theorem}{Theorem}

\newtheorem{example}{Example}

\catcode`\@=11
\def\downparenfill{$\m@th\braceld\leaders\vrule\hfill\bracerd$}
\def\overparen#1{\mathop{\vbox{\ialign{##\crcr\crcr \noalign{\kern0.4ex}
\downparenfill\crcr\noalign{\kern0.4ex\nointerlineskip}
$\hfil\displaystyle{#1}\hfil$\crcr}}}\limits}
\catcode`\@=12

\makeatletter
\newcommand{\pushright}[1]{\ifmeasuring@#1\else\omit\hfill$\displaystyle#1$\fi\ignorespaces}
\newcommand{\pushleft}[1]{\ifmeasuring@#1\else\omit$\displaystyle#1$\hfill\fi\ignorespaces}
\makeatother

\def\RR{{\mathbb R}}

\def\gz{\mathfrak{g}}

\def\Ouvs{\mathcal{S}} 
\def\B{\mathcal{B}}
\def\C{\mathcal{C}}

\def\R{\mathcal{R}}
\def\U{\mathcal{U}}

\def\blacksquare{\vrule height 0.5em depth 0pt width 0.5em}

\def\A{\mathcal{A}}
\def\B{\mathcal{B}}
\def\P{\mathcal{P}}

\def\H{\mathbf{H}}

\def\dz{{d_z}}
\def\du{{d_u}}  
\def\dx{{d_x}} 
\def\od{{o_d}} 
\def\ot{{o_t}}  

\DeclareMathAlphabet{\mathcursive}{U}{esstixcal}{m}{n}
\DeclareMathAlphabet\EuScript{U}{eus}{m}{n}
\SetMathAlphabet\EuScript{bold}{U}{eus}{b}{n}
\def\od{{\mathchoice%
{{\mbox{$\scriptstyle \EuScript{O}$}}}%
{{\mbox{$\scriptstyle \EuScript{O}$}}}%
{{\mbox{$\scriptscriptstyle \EuScript{O}$}}}%
{{\mbox{$\scriptscriptstyle \EuScript{O}$}}}%
}}
\def\ot{{\mathchoice%
{{\mbox{$\scriptstyle \EuScript{T}$}}}%
{{\mbox{$\scriptstyle \EuScript{T}$}}}%
{{\mbox{$\scriptscriptstyle \EuScript{T}$}}}%
{{\mbox{$\scriptscriptstyle \EuScript{T}$}}}%
}}

\definecolor{LPcomcolor}{rgb}{1,0,0}
\definecolor{PBcomcolor}{rgb}{1,0,1}
\definecolor{VAcomcolor}{rgb}{0,0,1}
\definecolor{removecolor}{cmyk}{1,0,0,0}

\long\def\supprimer#1{}

\begin{document}

\begin{frontmatter}

\title{On the triangular canonical form for uniformly observable controlled systems} 

\author[First]{P. Bernard}\ead{pauline.bernard@mines-paristech.fr},   
\author[First]{L. Praly}\ead{laurent.praly@mines-paristech.fr},     
\author[Second]{V. Andrieu}\ead{vincent.andrieu@gmail.fr},   
\author[Second]{H. Hammouri}\ead{hassan.hammouri@univ-lyon1.fr}  

\address[First]{Centre Automatique et Syst\`{e}mes, MINES ParisTech, PSL Research University, France }  
\address[Second]{Universit\'e Lyon 1, Villeurbanne, France -- CNRS, UMR 5007, LAGEP, France}


\begin{abstract}                        
We study controlled systems which are uniformly observable
and differentially observable with an order larger than the
system state dimension.
We establish that they may be transformed into a (partial) triangular canonical form but with possibly non locally Lipschitz functions. We characterize the points where this Lipschitzness may be lost and investigate the link with uniform infinitesimal observability.
\end{abstract}

\end{frontmatter}

\section{Introduction}

\subsection{Context}

A lot of attention has been dedicated to the construction of nonlinear
observers.  Although a general theory has been obtained for linear
systems, very few general approaches exist for nonlinear systems.  In
particular, the theory of high gain (\cite{KhaPra} and references therein) and Luenberger
(\cite{AndPra,And}) observers have been developed for autonomous non
linear systems but their extension to controlled systems is not
straightforward.

For designing an observer for a system, a preliminary step is often 
required. It consists in finding a reversible coordinate 
transformation, allowing us to rewrite the system dynamics in a 
targeted form 
more favorable for writing and/or analyzing the observer. In 
presence of input, two tracks are possible depending on whether we 
consider the input as a simple time function, making the system time 
dependent or as a more involved infinite dimensional parameter, making the system a family of dynamical 
systems, indexed by the input. In the former case, the transformation 
mentioned above is considered time dependent, and thus may need to be redesigned for each input. In the 
latter case, the transformation can be input-dependent.
Specifically :
\begin{list}{}{%
\parskip 0pt plus 0pt minus 0pt%
\topsep 0.5ex plus 0pt minus 0pt%
\parsep 0pt plus 0pt minus 0pt%
\partopsep 0pt plus 0pt minus 0pt%
\itemsep 0pt plus 0pt minus 0pt
\settowidth{\labelwidth}{-}%
\setlength{\labelsep}{0.5em}%
\setlength{\leftmargin}{\labelwidth}%
\addtolength{\leftmargin}{\labelsep}%
}
	\item[-] in  \cite{JouGau} (see also \cite{GauKup}), the transformation depends on the inputs and its derivatives. When the \textit{ACP(N) condition} is verified (see Lemma \cite[Definition 5.2.1, Lemma 5.2.2]{GauKup}), it leads to the so called phase-variable representation as targeted form (see \cite[Definition 2.3.1]{GauKup}), for which a high gain observer can be built.
	\item[-] in \cite{Bes}, the transformation does not depend on the input, and leads to a block triangular form when the system verifies the observability rank condition (see \cite{HerKre}). However, afterwards, an extra condition on the input is needed for the observer design.
	\item[-] in \cite{GauBor,GauHamOth}, the transformation does not depend on the input, and leads to a 
	triangular form when the system is a) uniformly observable (see \cite[Definition I.2.1.2]{GauKup} or 
Definition \ref{uniObs} below), and b) strongly differentially observable of order equal to the system state dimension
(see Definition \ref{diffObs} 
below). This so-called observable canonical form allows the design of a high gain observer.
\end{list}

In this paper, we complete and detail the results announced in \cite{BerPraAndNOLCOS}. We work within the third context
(of the second track),
but going beyond \cite{GauKup} with allowing the strong
differential observability order to be larger than the 
system state dimension.
We shall see that, in this case, the system dynamics may still be 
described by a
(partial) triangular canonical form (see (\ref{3}))
but with functions which may 
be non locally Lipschitz. 

\subsection{Definitions and problem statement}

We consider a controlled system of the form :
\begin{equation}
\dot{x}=f(x)+g(x) u \quad,\qquad y=h(x) 
\label{eqSyst}
\end{equation}
where $x$ is the state in $\RR^\dx $, $u$ is an input in $\RR^{\du }$,   $y$ 
is a measured output in $\RR$
and the functions $f$, $g$ and $h$ are sufficiently many times differentiable, $f$ being a column $\dx$-dimensional vector field
and $g$ a $(\dx \times \du )$-dimensional matrix field.
In the following, for a scalar function $\alpha $, $L_f\alpha $ denotes its Lie derivative in the direction of 
$f$. It has scalar values. We denote 
\begin{equation}
\label{LP1}
\textbf{H}_i(x) = (h(x),L_fh(x),...,L_f^{i-1}h(x)) 
\in \RR^i\ .
\end{equation}
It is a column $i$-dimensional vector. Similarly $L_g\alpha $ denotes the Lie derivative along each of the $\du $ 
columns of $g$. It has row $\du$-dimensional vector values.

Given an input time function $t\mapsto u(t)$ taking values in a compact subset $\U$ of $\RR^\du $, we
denote $X_u(x,t)$ a solution of (\ref{eqSyst}) going
through $x$ at time $0$.
We are
interested in solving :

\textbf{Problem $\mathbf{\P}$ : }
\textit{
Given a compact subset $\C$ of $\RR^\dx $, under which condition do there exist integers $\ot$ and $\dz$, a continuous injective function
$\Psi:\C\to \RR^\dz $, and continuous functions $\varphi_{\dz}:\RR^{\dz}\to \RR$ and $\gz_i:\RR^{i}(\textrm{or }\RR^{\dz})\to \RR^{\du}$ such that, when $x$ is in $\C$ and satisfies  (\ref{eqSyst}) and $u$ is in $\U$,
$z=\Psi(x)$ satisfies
\begin{equation}
\label{3}
\begin{array}{rcl}
\dot{z}_1&=&z_2 + \gz_1(z_1) \, u
\\[-0.5em]&\vdots&  \\[-0.5em]
\dot{z}_\ot &=&z_{\ot +1} + \gz_\ot (z_1,\ldots,z_\ot ) \, u  \\
\dot{z}_{\ot +1}&=&z_{\ot +2} + \gz_{\ot +1}(z) \, u  
\\[-0.5em]&\vdots&  \\[-0.5em]
\dot{z}_\dz &=&\varphi_\dz (z) + \gz_\dz(z) \, u 
\end{array}
\quad , \quad y =x_1
\end{equation}
}

Because
$\gz_i$ depends only on $z_1$ to $z_i$, for $i\leq \ot$, but potentially on all the components of $z$ for $i> \ot$, 
we call this particular form up-to-$\ot $-triangular canonical form
and $\ot$ is called the order of triangularity. When $\dz=\ot +1 $, we say full triangular 
canonical form. When the functions $\varphi _\dz $ and $\gz_j$ are locally 
Lipschitz we say Lipschitz up-to-$\ot$-triangular canonical form.

We are interested in addressing the
Problem $\mathbf{\P}$ because, when 
the functions are Lipschitz and $\dz=\ot +1 $, we get the nominal form for which high gain observers can be designed and therefore
$X_u(x,t)$ can be estimated knowing $y$ and $u$
as long as $(X_u(x,t),u(t))$ is in the given compact set $\C\times \U$.

We will use the following two notions of observability defined in \cite{GauKup}~:
\begin{definition}[Differential observability\footnotemark]
\label{diffObs}
System (\ref{eqSyst}) is weakly differentially observable of order $\od $ on an open subset $\Ouvs$ of $\RR^\dx$ if the function
$\H_\od  $ (see (\ref{LP1}))
is injective on $\Ouvs$. If it is also an immersion, the system is called strongly differentially observable of order $\od $. 
\end{definition}
\footnotetext{This notion is weaker than the usual differential observability defined for instance in \cite[Definition  I.2.4.2]{GauKup} for controlled systems, because it is a differential observability of the drift system only, namely when $u\equiv0$.}

\begin{definition}[Uniform observability]
\hskip -1.3pt (See \cite[Definition  I.2.1.2]{GauKup}.)
\label{uniObs}
System (\ref{eqSyst}) is uniformly observable on an open subset $\Ouvs$ of $\RR^\dx$ if, for any pair $(x_a,x_b)$ in $\Ouvs^2$ with $x_a\neq x_b$,  any strictly positive number $T$, and any $C^1$ function 
$u$ defined on $[0,T)$, there exists a time $t<T$ such that
$
h(X_u(x_a,t))\neq h(X_u(x_b,t)) 
$
and 
$
(X_u(x_a,s), X_u(x_b,s))\in \Ouvs^2 
$ 
for all $s \leq t$.
\end{definition}

Note that this notion is a matter of instantaneous observability since $T$ can be arbitrarily small. 
In the case
where $\H_\dx $ is a diffeomorphism,
we have

\begin{proposition}(See {\cite{GauBor,GauHamOth}})\label{Prop_Gau}
If System (\ref{eqSyst}) is uniformly observable and strongly differentially observable of order $\od =\dx $
on an open set $\Ouvs$ containing the given compact set $\C$,
it can be transformed on $\C$ into
a full Lipschitz triangular canonical form of dimension $\dz=\dx $.
\end{proposition}

In general, it is possible for the system not to be strongly differentially observable of order
$\dx$ everywhere.  This motivates our interest in the case where the system is strongly differentially observable of
order $\od >\dx $, i-e $\textbf H_\od $ is an injective immersion and not a diffeomorphism.   As
we shall see in Section \ref{sec_immCase}, in this case, we may still have an at least
up-to-$(\dx+1) $-triangular form, but the Lipschitzness of its nonlinearities can be lost.  Since
this property is crucial for the implementation of high gain observers (see \cite{CicDalGer}),
we give in Section \ref{secLip} some sufficient conditions under which the
Lipschitzness is preserved.


\section{Immersion case ($\od >\dx $)}
\label{sec_immCase}

The specificity of the triangular canonical form
 (\ref{3})
 is not so much in 
its structure but more in the dependence of its functions $\gz_i$ and 
$\varphi _\dz $. Indeed, by choosing $\Psi=\textbf{H}_{\dz}$, 
 we obtain:
\\[0.7em]\vbox{\noindent$\displaystyle 
\dot{\overparen{\H_{\dz}(x)}}=
$\hfill \null \\\null \hfill $\displaystyle \left(\begin{array}{ccccc}
0 & 1 &0&\ldots &0
\\[-0.5em]
\vdots & \ddots& \ddots&\ddots&\vdots
\\[-0.5em]
\vdots & & \ddots&\ddots&0
\\
0 & \ldots & \ldots&0 & 1
\\
0 & \ldots &\ldots & \ldots& 0
\end{array}\right)\H_{\dz}(x)
+
\left(\begin{array}{c}
0 \\[-0.5em] \vdots\\[-0.5em]\vdots\\[-0.5em]0\\L_f^{\dz}h(x) 
\end{array}\right)
+
L_g\H_\dz(x) u
$}\\[0.7em]
To get (\ref{3}),
we need further the existence of a 
function $\varphi _\dz$ satisfying
\\[0.5em]\null \hfill $\displaystyle 
L_f^{\dz}h(x)\;=\; \varphi _\dz(\H_\dz(x)) \quad \forall x\in \C
$\refstepcounter{equation}\label{LP16}\hfill$(\theequation)$\\[0.5em]
and, for $i\leq \ot$, of  functions $\gz_i$ satisfying
\\[0.5em]\null \hfill $\displaystyle 
L_gL_f^{i-1}h(x)\;=\; \gz_i(h(x),\ldots,L_f^{i-1}h(x))  \quad \forall x\in \C \ .
$\refstepcounter{equation}\label{LP17}\hfill$(\theequation)$\\[0.5em]
Let us illustrate via the following elementary example what can occur. 
\begin{example}
\label{ex1}
Consider the system defined as
\\[0.5em]\null \hfill $\displaystyle 
\dot{x}_1 \;=\;  x_2 \  ,\quad 
\dot{x}_2 \;=\; x_3^3 \  ,\quad 
\dot{x}_3 \;=\; 1+u \  ,\quad 
y \;=\;  x_1
$\hfill \null \\[0.5em]
We get
\\[0.5em]\null \hfill $ 
\begin{array}{r@{\; }c@{\; }l}
\H_3(x)&=&(h(x),L_fh(x),L_f^2h(x)) = (x_1,x_2,x_3^3)
\\
\H_5(x)&=&(\H_3(x),L_f^3h(x),L_f^4h(x)) 
=(\H_3(x),3x_3^2,6x_3)
\end{array}
$\hfill \null \\[0.5em]
Hence $\H_3$ is a bijection and $\H_5$
is an injective immersion on $\RR^3$. So this system is
weakly differentially observable of order $3$ on $\RR^3$
and strongly differentially observable of order $5$ on $\RR^3$.
Also the function $(x_1,x_2,x_3)\mapsto (y,\dot y,\ddot y)$ being injective for all $u$,
it is uniformly observable on $\RR^3$.
From this we could be tempted to pick $\dz=3$ or $5$ and the compact set $\C$ arbitrary in $\RR^3$.
Unfortunately, if we choose $\dz=3$, we have
\\[0.5em]\null \hfill $\displaystyle 
L_f^{3}h(x)\;=\; 3x_3^2\;=\; 3(L_f^2h(x))^{2/3}
$\hfill \null \\[0.5em]
and there is no locally
Lipschitz function $\varphi _3$ satisfying 
(\ref{LP16}) if the given compact set $\C$ contains a point satisfying $x_3=0$.
If we choose $\dz=5$, we have
\\[0.5em]\null \hfill $\displaystyle 
L_gL_f^2h(x)=3x_3^2=L_f^3h(x)=3(L_f^2h(x))^{2/3}
$\hfill \null \\[0.5em]
and there is no locally Lipschitz function $\gz _3$ satisfying 
(\ref{LP17})  if the given compact set $\C$ contains a point satisfying $x_3=0$. 
\end{example}

Leaving aside the Lipschitzness requirement for
the time being,  we focus  
on the existence of continuous functions $\varphi _\dz $ and $\gz_i$ verifying (\ref{LP16}) and (\ref{LP17}) respectively. 

\subsection{Existence of $\varphi _\dz $ satisfying (\ref{LP16})}
\begin{proposition}
	\label{propvarphi}
	Suppose System (\ref{eqSyst}) is weakly differentially observable of order 
	$\od $ on an open set $\Ouvs$ containing the given
	compact set $\C$. For any $\dz\geq \od$, there 
	exists a continuous function $\varphi _\dz :\RR^\dz \to \RR$ satisfying (\ref{LP16}). If System (\ref{eqSyst}) is strongly differentially observable of order 
	$\od$ on $\Ouvs$, the  function
	$\varphi _\dz $ can be chosen Lipschitz on $\RR^\dz $.
\end{proposition}
\begin{pf}
There is nothing really new in this result. It is a direct consequence of the fact that a continuous injective function defined on a 
compact set admits a continuous left inverse defined on the image (see  \cite[\S 16.9]{Bartle}) and that this 
left inverse can be extended to the full space (e.g. Tietze extension theorem). 
In the case where $\H_{\dz}$ is also an immersion,
according to \cite[Lemma 3.2]{And}, there exists a real number $L_H>0$ such that
\\[0.5em]\null \hfill $\displaystyle 
|x_a-  x_b|\leq  L_H\left|\textbf{H}_\dz (x_a)-\textbf{H}_\dz (x_b)\right|
\quad \forall (x_a,x_b)\in \C^2\; .
$\hfill \null \\[0.5em]
Therefore, the previously mentioned continuous left-inverse of $\H_{\dz}$ defined on $\H_{\dz}(\C)$ is Lipschitz on $\H_{\dz}(\C)$. According to \cite[Theorem 1]{McS}, it admits a Lipschitz extension defined on $\RR^\dz$.
\hfill $\blacksquare$
\end{pf}

\subsection{Existence of $\gz_i$ satisfying (\ref{LP17})}

Concerning the functions $\gz_i$, we will prove the following result :

\begin{proposition}
	\label{claimGi}
	Suppose System (\ref{eqSyst}) is uniformly observable on an open set $\Ouvs$ containing the given compact set
	$\C$.
	\begin{list}{}{%
			\parskip 0pt plus 0pt minus 0pt%
			\topsep 0.5ex plus 0pt minus 0pt%
			\parsep 0pt plus 0pt minus 0pt%
			\partopsep 0pt plus 0pt minus 0pt%
			\itemsep 0.5ex plus 0pt minus 0pt
			\settowidth{\labelwidth}{-}%
			\setlength{\labelsep}{0.5em}%
			\setlength{\leftmargin}{\labelwidth}%
			\addtolength{\leftmargin}{\labelsep}%
		}
		\item[-] There exists a continuous function
		$\gz_1:\RR\to \RR^{\du }$ satisfying (\ref{LP17}).
		\item[-] If, for some $i$ in $\{2,\ldots,\dx\}$,  $\H_2, \dots, \H_i$ defined in \eqref{LP1} are open maps, then,  for all $j \leq i$, there exists a continuous function $\gz_j:\RR^j\to \RR^{\du }$ satisfying (\ref{LP17}).
	\end{list} 
\end{proposition}

The rest of this section is dedicated to the proof of this result through a series of lemmas, the proof of which can be found in appendix.

A first important thing to notice is that the following property must be satisfied 
for the identity \eqref{LP17} to be satisfied (on $\Ouvs$).

\textbf{Property $\mathbf{\A(i)}$:} 
\\[0.5em]$\displaystyle
L_gL_f^{i-1}h(x_a)=  L_gL_f^{i-1}h(x_b)
$\hfill \null
\\[0.3em]\null\hfill$\displaystyle
\forall (x_a,x_b)\in\Ouvs^2  \, : \, \textbf{H}_i(x_a)=\textbf{H}_i(x_b)
$

Actually the converse is true and is a direct consequence from Lemma \ref{lemTech2} proved in Appendix \ref{appTechLem} :

\begin{lem}
\label{lem1}
If Property $\A(i)$
is satisfied with $\Ouvs$ containing the given compact set $\C$, then 
there exists a continuous function $\gz_i:\RR^i\to \RR^{\du }$ satisfying (\ref{LP17}).
\end{lem}

Property $\A(i)$ being sufficient to obtain
the existence of a function $\gz_i$ satisfying (\ref{LP17}), we study now under which conditions it holds.
Clearly $\A(i)$ is
satisfied for all $i\geq \od$ if $\textbf{H}_\od$ is injective.  If we do not have this injectivity property 
the situation is more complex. To overcome the difficulty we introduce the following assumption for $2\leq i\leq d_x+1$.

\textbf{Property $\mathbf{\B(i)}$ : }
\refstepcounter{equation}\label{prop_Bi}\hfill$(\theequation)$
\\[0.2em]
\textit{%
For any $(x_a,x_b)$ in $\Ouvs^2$ such that $x_a\neq x_b$, verifying
$
\textbf{H}_i(x_a)=\textbf{H}_i(x_b) \ ,
$	
there exists a sequence $(x_{a,k},x_{b,k})$ of points in $\Ouvs^2$ converging to $(x_a,x_b)$ such that for all $k$,
$\textbf{H}_i(x_{a,k})=\textbf{H}_i(x_{b,k})$
and $\frac{\partial \textbf{H}_{i-1}}{\partial x}$ is full-rank at $x_{a,k}$ or $x_{b,k}$.
}

As in this Property, let $x_a\neq x_b$ be such that $\textbf{H}_i(x_a)=\textbf{H}_i(x_b)$.  If $\frac{\partial
\textbf{H}_{i-1}}{\partial x}$ is full-rank at either $x_{a}$ or $x_{b}$, then we can take $(x_{a,k},x_{b,k})$
constant equal to $(x_a,x_b)$.  Thus, it is sufficient to check $\B(i)$ around points where neither
$\frac{\partial \textbf{H}_{i-1}}{\partial x}(x_a)$ nor $\frac{\partial \textbf{H}_{i-1}}{\partial x}(x_b)$ is
full-rank.  But according to \cite[Theorem 4.1]{GauKup}, the set of points where $\frac{\partial
\textbf{H}_{\dx}}{\partial x}$ is not full-rank is of codimension at least one for a uniformly observable
system.  Thus, it is always possible to find points $x_{a,k}$ as close to $x_a$ as we want such that
$\frac{\partial \textbf{H}_{i-1}}{\partial x}(x_{a,k})$ is full-rank.  The difficulty of $\B(i)$ thus rather
lies in ensuring that we have also $\textbf{H}_i(x_{a,k})=\textbf{H}_i(x_{b,k})$. 

In Appendix \ref{appPrLemclaimUniObs}, we prove
\begin{lem}
\label{claimUniObs}
Suppose System (\ref{eqSyst}) is uniformly observable on a set $\Ouvs$. 
\begin{list}{}{%
		\parskip 0pt plus 0pt minus 0pt%
		\topsep 0.5ex plus 0pt minus 0pt%
		\parsep 0pt plus 0pt minus 0pt%
		\partopsep 0pt plus 0pt minus 0pt%
		\itemsep 0.5ex plus 0pt minus 0pt
		\settowidth{\labelwidth}{-}%
		\setlength{\labelsep}{0.5em}%
		\setlength{\leftmargin}{\labelwidth}%
		\addtolength{\leftmargin}{\labelsep}%
	}
	\item[-] Property $\A(1)$ is satisfied.
	\item[-] If, for some $i$ in $\{2,\ldots,\dx +1\}$, Property
	$\B(i)$  holds and Property $\A(j)$ is satisfied
for all $j$ in $\{1,\ldots,i-1\}$, then Property $\A(i)$ holds.
\end{list} 
\end{lem}

Thus, the first point in Proposition \ref{claimGi} is proved. Besides, a direct consequence of Lemmas \ref{lem1} and \ref{claimUniObs} is that a sufficient condition to have the existence of the functions $\gz_i$ for $i$ in $\{2,\ldots,\dx+1\}$ is to have $\B(j)$ for $j$ in $\{2,\ldots,i\}$.  The following lemma finishes the proof of Proposition \ref{claimGi} by showing that $\B(j)$ is in fact satisfied when $\H_j$ is an open map.

\begin{lem}
Suppose that for some $j$ in $\{2,\ldots,\dx\}$, $\H_j$	is an open map on $\Ouvs$. Then, $\B(j)$ is satisfied.
\end{lem}

\begin{pf}
Take $(x_a,x_b)$ in $\Ouvs^2$ such that $x_a\neq x_b$ and $\textbf{H}_j(x_a)=\textbf{H}_j(x_b)=y_0$. Let $\Pi$ be the set of points of $\Ouvs$ such that $\frac{\partial \textbf{H}_{j}}{\partial x}$ is not full-rank. According to Sard's theorem, $\H_j(\Pi)$ is of measure zero in $\RR^j$. Now, take $p>0$ and consider $B_p(x_a)$ and $B_p(x_b)$ the open balls of radius $\frac{1}{p}$ centered at $x_a$ and $x_b$ respectively. Since $\H_j$ is open, $\H_j(B_p(x_a))$ and $\H_j(B_p(x_ b))$ are open sets, both containing $y_0$. Thus,  $\H_j(B_p(x_a))\cap\H_j(B_p(x_ b))$ is a non-empty open set. It follows that $(\H_j(B_p(x_a))\cap\H_j(B_p(x_ b)))\setminus \H_j(\Pi)$ is non-empty and contains a point $y_p$. We conclude that there exist $(x_{a,p},x_{b,p})$ in $B_p(x_a)\times B_p(x_b)$ such that $\textbf{H}_j(x_{a,p})=\textbf{H}_j(x_{b,p})=y_p$ and $\frac{\partial \textbf{H}_{j}}{\partial x}$ (and thus $\frac{\partial \textbf{H}_{j-1}}{\partial x}$) is full-rank at $x_{a,p}$ and $x_{b,p}$.  Besides $(x_{a,p},x_{b,p})$ converges to $(x_a,x_b)$, and $\B(j)$ is satisfied.
\hfill $\blacksquare$
\end{pf}

Note that the assumption that $\H_j$ is an open map is stronger that $\B(j)$ since it leads to the full rank of $\frac{\partial \textbf{H}_{j}}{\partial x}$, while, in $\B(j)$, we only need the full-rank for $\frac{\partial \textbf{H}_{j-1}}{\partial x}$. We show in the following example that the openness of $\H_j$ is not necessary.

\begin{example}
	\label{ex2}
Consider the system defined as
\\[0.5em]\null \hfill $\displaystyle 
\dot{x}_1 \;=\;  x_2 \  ,\quad 
\dot{x}_2 \;=\; x_3^3x_1 \  ,\quad 
\dot{x}_3 \;=\; 1+u \ , \quad
y \;=\; x_1
\refstepcounter{equation}(\arabic{equation}) \label{eq_dyn_ex2}$\hfill \null \\[0.5em]
On $\Ouvs= \left\{x\in\RR^3:\,  x_1^2+x_2^2 \neq 0\right\}$, and whatever $u$ is, the knowledge of the function $t\mapsto y(t) = X_1(x,t)$ and therefore of its three first derivatives
\\[0.5em]\null \hfill $\displaystyle
\dot{y} \;=\; x_2 \ , \quad
\ddot{y} \;=\; x_3^3x_1 \ , \quad
\dddot{y}\;=\; 3x_3^2	x_1(1+u) + x_3^3 x_2
$\hfill \null \\[0.5em]
gives us $x_1$, $x_2$ and $x_3$. Thus, the system is uniformly observable on $\Ouvs$. 
Besides, the function
\\[0.5em]\null \hfill $\displaystyle
\mathbf{H}_4(x) = \left(
x_1 \, , \, x_2 \, , \, x_3^3x_1 \, , \, 3x_3^2x_1 + x_3^3x_2 
\right)
$\hfill \null \\[0.5em]
is injective on $\Ouvs$ and the system is weakly differentially observable of order 4 on $\Ouvs$. Now, although $\H_2$ is trivially an open map on $\Ouvs$, $\H_3$ is not. Indeed, consider for instance the open ball $B$ in $\RR^3$ with radius $\frac{1}{2}$ and centered at $(0,x_2,0)$   for some $x_2$ such that $|x_2|> \frac{1}{2}$. $B$ is contained in $\Ouvs$. Suppose its image by $\H_3$ is an open set of $\RR^3$. It contains $\H_3(0,x_2,0)=(0,x_2,0)$ and thus $(\varepsilon,x_2,\varepsilon)$ for any sufficiently small $\varepsilon$. This means that there exist $x$ in $B$ such that $(\varepsilon,x_2,\varepsilon)= \H_3(x)$, i-e necessarily $x_1=\varepsilon$ and $x_3=1$. But this point is not in $B$, and we have a contradiction. Therefore, $\H_3$ is not open.
However, $\B(3)$ trivially holds because $\H_2$ is full-rank everywhere. 
\end{example}

\subsection{A solution to Problem $\P$}

With Propositions \ref{propvarphi} and \ref{claimGi}, we have the following solution to Problem $\P$.
\begin{theorem}
	\label{thm}
Suppose System (\ref{eqSyst}) is weakly differentially observable  of order $\od$ and uniformly observable  
on an open set $\Ouvs$ containing the given compact set $\C$. With selecting $\Psi=\mathbf{H_\od}$ and $\dz = \od$,
we have a solution to Problem $\P$ if we pick either $\ot = 1$, or
$\ot = i$ when $\H_j$ is an open map for any $j$ in $\{2,\dots,i\}$ with $i\leq \dx$.
\end{theorem}
\begin{pf}
In each case, the function $\varphi_\dz$ is obtained from Proposition \ref{propvarphi}. Functions
$\gz_i$ for $i\leq \ot$ are obtained from Proposition \ref{claimGi}. Finally, it is  possible to construct the 
functions $\gz_i$, for  $i > \ot$, in the same way as $\varphi_\dz$ in Proposition \ref{propvarphi}.
\end{pf}

Three remarks :
\begin{list}{}{%
		\parskip 0pt plus 0pt minus 0pt%
		\topsep 0.4ex plus 0pt minus 0pt%
		\parsep 0pt plus 0pt minus 0pt%
		\partopsep 0pt plus 0pt minus 0pt%
		\itemsep 0.5ex plus 0pt minus 0pt
		\settowidth{\labelwidth}{-}%
		\setlength{\labelsep}{0.5em}%
		\setlength{\leftmargin}{\labelwidth}%
		\addtolength{\leftmargin}{\labelsep}%
} 
\item[-] As seen in Example \ref{ex2}, the openness of the functions $\H_j$ is sufficient but not necessary. We may ask only for $\B(j)$ for any $j$ in $\{2,\dots,i\}$ with $i\leq \dx+1$. Besides, this weaker assumption allows to obtain the existence of $\gz_i$ up to the order $\dx+1$.
\item[-] Consider the case where $\B(j)$ is satisfied for all $j\leq \dx+1$ and $\od=\dx+2$.  Then
we have $\ot=\dx+1$ and it is possible to obtain a full triangular form of dimension $\dz=\ot+1=\od=\dx+2$.
Actually, we still have a full triangular form if we choose $\dz>\od$.  Indeed, $\textbf{H}_\od$
being injective, $\A(i)$ is satisfied for all $i$ larger than $\od$, thus there also exist continuous
functions $\gz_i:\RR^i\to \RR^{\du }$ satisfying (\ref{LP17}) for all $i\geq \od$.  It follows that $\ot$ can be taken larger than $\dx+1$ and $\dz=\ot+1$ larger than $\od$.
\item[-] If $\P$ is solved with $\dz=\ot+1$, we have a full triangular canonical form of
dimension $\dz$.  But, at this point we know nothing about the regularity of the functions $\gz_i$,
besides continuity.  As we saw in Example \ref{ex1}, even the usual assumption of strong differential
observability is not sufficient to make it Lipschitz everywhere.  This may impede (see for instance
\cite{CicDalGer}) the convergence of the high gain observer proposed in \cite{GauHamOth}.  That is
why, in the next section, we look for conditions under which the Lipschitzness is ensured. 
\item[-]
As mentioned in the introduction, another way of solving Problem $\P$ is to allow the transformation
$\Psi$ to depend on the control $u$ and its derivatives.  In particular, if $\dz > \ot+1$, a full
triangular form may still be obtained with $\Psi = (\textbf{H}_{\ot}, \tilde{\Psi})$ where the
components $\tilde{\Psi}_i$ of $\tilde{\Psi}$ are defined recursively as
\\[0.4em]\null \hfill $\displaystyle 
\tilde{\Psi}_1 = L_f^\ot h \quad , \quad
\tilde{\Psi}_{i+1} = L_{f+gu}\tilde{\Psi}_i +\sum_{j=0}^{i-2} \frac{\partial \tilde{\Psi}_i}{\partial u^{(j)}}u^{(j+1)} 
$\hfill \null \\[0.5em]
until (if possible) the map $x\mapsto \Psi(x,u,\dot{u},...)$ becomes injective for all $(u,\dot{u},...)$.
\end{list}

\begin{example}
	\label{ex2_suite1}
	Coming back to Example \ref{ex2}, we have seen that $\H_2$ is open and that $\H_3$ is not but $\B(3)$ is satisfied. Besides, the system is weakly differentially observable of order 4. We deduce that there exists a full-triangular form of order 4. Indeed, we have $L_gh(x) = L_gL_fh(x) = 0$ and
	\\[0.5em]\null \hfill $\displaystyle 
	L_gL_f^2h(x) = 3x_3^2x_1 = 3 (L_f^2h(x))^{\frac{2}{3}} (h(x))^{\frac{1}{3}} 
	$\hfill \null \\[0.5em]
	so that we can take 
	\\[0.5em]\null \hfill $\displaystyle 
	\gz_1=\gz_2=0 \quad , \quad \gz_3(z_1,z_2,z_3) = 3 z_3^{\frac{2}{3}} z_1^{\frac{1}{3}}  . 
	$\hfill \null \\[0.5em]
	As for $\varphi_4$ and $\gz_4$, they are obtained via inversion of $\H_4$   i-e for instance on $\RR^4\setminus\{(0,0,z_3), \, z_3\in \RR\}$
	\\[0.5em]\null \hfill $\displaystyle
	\H_4^{-1}(z) = \left(z_1,z_2,\left(\frac{(z_4-3 z_3^{\frac{2}{3}} z_1^{\frac{1}{3}} )^2+z_3^2}{z_1^2+z_2^2} \right)^{\frac{1}{6}}\right) \ .
	$\hfill \null \\[0.5em]
\end{example}


\section{Lipschitzness of the triangular form}
\label{secLip}

\subsection{A sufficient condition}
We saw with Examples \ref{ex1} and \ref{ex2} that uniform observability is not sufficient
for the functions $\gz_i$ to be Lipschitz. Nevertheless,
we are going to see in this section that it is sufficient except maybe around the image of points where $\displaystyle \frac{\partial\H_{i}}{\partial x}$ is not full-rank ($x_1=0$ or $x_3=0$ in Example \ref{ex2}).

Consider the set $\R_i$ of points in $\Ouvs$ where $\displaystyle\frac{\partial\H_i}{\partial x}$ has
full rank.  Note that according to \cite[Corollaire p68-69]{Leb}, if $\H_i$ is an open map, $\R_i$ is an open dense set.
Assume $\R_\dx\cap \C$ is non empty. Then there exists $\varepsilon_0>0$ such that, for all $\varepsilon$ in $(0,\varepsilon_0]$, the set
\\[0.5em]\null Ê\hfill $\displaystyle 
K_{i,\varepsilon } = \left\{x\in \R_i\cap \C \ ,  \quad d(x\,  ,\,  
\RR^\dx \!\setminus{\!\R_i}
) \geq \varepsilon \right\} \ .
$\hfill \null \\[0.5em]
is non-empty and compact, and such that its points are ($\varepsilon$)-away from singular points. The next proposition shows that the functions $\gz_i$ can be taken Lipschitz on the image of $K_{i,\varepsilon }$, i-e everywhere except arbitrary close to the image of points where the rank 
of the Jacobian of $\H_i$ drops.  

\begin{proposition}
	\label{propRgCst_glob}
	Assume System (\ref{eqSyst}) is uniformly observable on an open set $\Ouvs$ containing the compact set $\C$. For all $i$ in $\{1,...,\dx\}$ and for any $\varepsilon$ in $(0,\varepsilon_0]$, there exists a Lipschitz function $\gz_i:\RR^i\to\RR^\du$ satisfying (\ref{LP17}) for all $x$ in $K_{i,\varepsilon}$.
\end{proposition}

\begin{pf}
	As noticed after the statement of Property $\B(i)$, since
	$\displaystyle\frac{\partial\H_i}{\partial x}$ has
	full rank in the open set $\R_i$, Property $\B(i)$ holds on $\R_i$ (i-e with $\R_i$ replacing $\Ouvs$ in its 
	statement). It follows from Lemma \ref{claimUniObs} that $\A(i)$ is satisfied on $\R_i$.
	Besides, according to Lemma \ref{lemRgCst}, $\H_i(\R_i)$ is open and there exists a $C^1$ function $\gz_i$ defined on $\H_i(\R_i)$ such that for all $x$ in $\R_i$, $\gz_i(\H_i(x))=L_gL_f^{i-1}h(x)$. Now, $K_{i,\varepsilon}$ being a compact set contained in $\R_i$, and $\H_i$ being continuous, $\H_i(K_{i,\varepsilon})$ is a compact set contained in $\H_i(\R_i)$. Thus, $\gz_i$ is Lipschitz on $\H_i(K_{i,\varepsilon})$. According to \cite{McS}, there exists a Lipschitz extension of $\gz_i$ to $\RR^i$ coinciding with $\gz_i$ on  $\H_i(K_{i,\varepsilon})$, and thus verifying (\ref{LP17}) for all $x$ in $K_{i,\varepsilon}$.
	\hfill $\blacksquare$
\end{pf}

For a strongly differentially observable system of order $\od =\dx$ on $\Ouvs$,
the Jacobian of $\H_i$ for any $i$ in $\{1,...,\dx\}$ has full rank on $\Ouvs$. Thus, taking
$\dz=\ot+1=\od=\dx$ a full Lipschitz triangular form of dimension $\dx$ exists, i.e.
we recover  the result of Proposition \ref{Prop_Gau}. 

\begin{example}
	\label{ex2_suite2}
	In Example \ref{ex2}, $\H_3$ is full rank on $\Ouvs\setminus \{ x\in \RR^3 \ | \ x_1=0 \, \textrm{or } x_3=0 \}$. Thus, according to Proposition \ref{propRgCst_glob}, the only points where $\gz_3$ may not be Lipschitz, are the image of points where $x_1=0$ or $x_3=0$. Let us study more precisely what happens around those points. Take $x_a=(x_{1,a},x_{2,a},0)$ in $\Ouvs$. If  there existed a locally Lipschitz function $\gz_3$ verifying \eqref{LP17} around $x_a$, there would exist $L>0$ such that for any $x_b=(x_{1,b},x_{2,a},x_{3,b})$ sufficiently close to $x_a$ with $x_{1,b}\neq 0$,   $|3x_{3,b}^2|\leq L |x_{3,b}^3|$, which we know is impossible. Therefore, there does not exist a function $\gz_3$ which is Lipschitz around the image of points where $x_3=0$. Let us now study what happens elsewhere, namely on $\tilde{\Ouvs}=\Ouvs \setminus \{ x\in \RR^3 \ |  x_3=0 \}$. It turns out that on any compact set $\C$ of $\tilde{\Ouvs}$, there exists\footnote{
		If $x_{1,a}$ and $x_{1,b}$  are both zero, the inequality is trivial. Suppose $|x_{1,a}|> |x_{1,b}|$ and denote	$\rho=\frac{x_{1,b}}{x_{1,a}}$. If $\rho<0$, we have directly
		$|x_{3,a}^2-\rho \, x_{3,b}^2|\leq \max\{x_{3,a}^2,x_{3,b}^2\}|1-\rho|$. If now $\rho>0$, $x_{3,a}^2-\rho \, x_{3,b}^2=\frac{(x_{3,a}^3-\rho^{\frac{3}{2}} \, x_{3,b}^3)(x_{3,a}+\sqrt{\rho}x_{3,b})}{x_{3,a}^2+\sqrt{\rho}x_{3,a}x_{3,b}+\rho x_{3,b}^2}$ and thus $|x_{3,a}^2-\rho \, x_{3,b}^2|\leq \frac{2\sqrt{2}}{\sqrt{x_{3,a}^2+\rho x_{3,a}^2}}|x_{3,a}^3-\rho^{\frac{2}{3}} \, x_{3,b}^3|$. Besides, $|x_{3,a}^3-\rho^{\frac{3}{2}} \, x_{3,b}^3|= |x_{3,a}^3-\rho x_{3,b}^3 + \rho(1- \sqrt{\rho}) \, x_{3,b}^3|\leq |x_{3,a}^3-\rho x_{3,b}^3|+\frac{\rho |x_{3,b}^3|}{1+\sqrt{\rho}}|1-\rho|$ which gives $L$ on compact sets.
	}
	$L$ such that we have for all $(x_a,x_b)$ in $\C^2$,
	\\[0.7em]$\displaystyle 
	|x_{3,a}^2x_{1,a}-x_{3,b}^2x_{1,b}| \leq L (|x_{1,a}-x_{1,b}|
	$\hfill \null \\\null \hfill $\displaystyle 
	+ |x_{3,a}^3x_{1,a}-x_{3,b}^3x_{1,b}|)
	$\\[0.7em]
	Therefore, the continuous function $\gz_3$ found earlier in Example \ref{ex2_suite1} such that $\gz_3(\H_3(x))=L_gL_f^2(x)=3x_3^2x_1$ on $\Ouvs$ (and thus on $\C$)  verifies in fact
	$$
	|\gz_3(z_a)-\gz_3(z_b)| \leq L |z_a-z_b|
	$$
	on $\H_3(\C)$ and can be extended to a Lipschitz function on $\RR^3$ according to \cite[Theorem 1]{McS}. We conclude that although $\H_3$ does not have a full-rank Jacobian everywhere on $\C$ (singularities at $x_1=0$), it is possible to find a Lipschitz function $\gz_3$ solution to our problem on this set. 
\end{example}

\subsection{A necessary condition}

We have just seen that the condition in Proposition \ref{propRgCst_glob} that the Jacobian of $\H_i$ be full-rank,  is sufficient but not necessary. In order to have locally Lipschitz functions $\gz_i$ satisfying \eqref{LP17}, there must exist for all $x$ a strictly positive number $L$ such that  for all $(x_a,x_b)$ in a neighborhood of $x$,
\\[0.5em]\null \hfill$\displaystyle
|L_gL_f^{i-1}h(x_a)-L_gL_f^{i-1}h(x_b)| \leq L \,  |\H_i(x_a)-\H_i(x_b)| \ .
\ \hfill\null \refstepcounter{equation}(\arabic{equation}) \label{eq_inequalityLips} $\\[0.5em]
We have the following necessary condition :
\begin{lem}
	\label{lem_necCondLips}
	Consider $x$ in $\Ouvs$ such that \eqref{eq_inequalityLips} is satisfied in a neighborhood of $x$. Then, for any non zero vector $v$ in $\RR^\dx$, and any $k$ in $\{1,\ldots,\du\}$, we have :
	\begin{equation}
	\label{eq_condLipsJac}
	\frac{\partial \H_i}{\partial x}(x)\,  v=0 \quad \Rightarrow \quad \frac{\partial L_{g_k}L_f^{i-1}h}{\partial x}(x)\,  v = 0 \ .
	\end{equation}
\end{lem}

\begin{pf}
	Assume there exists a non-zero vector $v$ in $\RR^\dx$ such that
	$
	\displaystyle
	\frac{\partial \H_i}{\partial x}(x)\,  v=0 \  
	$. Choose $r>0$ such that Inequality \eqref{eq_inequalityLips} holds on $B_r(x)$, the ball centered at $x$ and of radius $r$.
	Consider for any integer $p$ the vector $x_p$ in $B_r(x)$ defined by
	$
	x_p=x-\frac{1}{p}\frac{r}{|v|}v .
	$
	This gives a sequence
	converging to $x$ when $p$ tends to infinity.
	We have
	\\[1em]$\displaystyle
	0\; \leq \; \frac{|L_{g_k}L_f^{i-1}h(x_p)-  L_{g_k}L_f^{i-1}h(x)|}{|x-x_p|}
	$\hfill \null
	\\[0.7em]\null\hfill$\displaystyle
	\leq \; L \;
	\frac{| \H_i(x_p)-\H_i(x)|}{|x-x_p|} \ 
	\refstepcounter{equation}(\arabic{equation}) \label{5} $\\[1em]
	The sequence $\frac{x-x_p}{|x-x_p|}$ tends to $v$, and thus, $\frac{\H_i(x_p)-\H_i(x)}{|x-x_p|}$ tends to $\frac{\partial \H_i }{\partial x}(x)\,  v$
	which by assumption is $0$
	.
	Similarly
	$\frac{1}{|x-x_p|}(L_{g_k}L_f^{i-1}h(x)-L_{g_k}L_f^{i-1}h(x_p))$ tends to $\displaystyle\frac{\partial L_{g_k}L_f^{i-1}h}{\partial x}(x)\,  v$
	which is also $0$ according to (\ref{5}).\hfill $\blacksquare$
\end{pf}

We conclude that when $\H_i$ does not have a full-rank Jacobian, it must satisfy condition \eqref{eq_condLipsJac} to allow the existence of locally Lipschitz triangular functions $\gz_i$. This condition is in fact about uniform infinitesimal observability.

\begin{definition}[Uniform infinitesimal observability]
	\hskip -1.3pt (See \cite[Definition  I.2.1.3]{GauKup}.)
	\label{uniInfObs}
	Consider the system lifted to the tangent bundle (see \cite[page 10]{GauKup})
	\begin{equation}
	\label{4}
	\left\{
	\begin{array}{rcl}
	\dot x &=& f(x) + g(x) u
	\\
	\dot{v}&=& \left[\frac{\partial f}{\partial x}(x) + \frac{\partial 
		gu}{\partial x}(x)\right] v
	\end{array}
	\right.
	\ , \quad
	\left\{
	\begin{array}{rcl}
	y &=& h(x)
	\\
	w&=& \frac{\partial h}{\partial x}(x) v
	\end{array}
	\right.
	\end{equation}
	with $v$ in $\RR^\dx $ and $w$ in $\RR$ and the solutions of which are
	denoted $(X_u(x,t),V_u((x,v),t))$.
	System (\ref{eqSyst}) 
	is uniformly infinitesimally observable on $\Ouvs$ if, for any pair $(x,v)$ in $\Ouvs\times\RR^\dx \setminus\{0\}$, any strictly positive number $T$, and any $C^1$ function $u$ defined on an interval $[0,T)$, there exists a time $t<T$ such that
	$\displaystyle
	\frac{\partial h}{\partial x}(X_u(x,t))\,  V_u((x,v),t)\neq 0 
	$
	and such that $X_u(x,s)\in \Ouvs$ for all $s\leq t$.
\end{definition}

We have the following result.
\begin{proposition}
	\label{prop_CN_infObs}
	Suppose that System (\ref{eqSyst}) is strongly differentially observable of order $\od$ (or at least that $\H_\od$ is an immersion on $\Ouvs$) and that Inequality \eqref{eq_inequalityLips} is verified at least locally around any point $x$ in $\Ouvs$ for any $i$ in $\{1,\ldots, \od \}$.
	Then the system is uniformly infinitesimally observable on $\Ouvs$.
\end{proposition}

\begin{pf}
	According to Lemma \ref{lem_necCondLips}, we have \eqref{eq_condLipsJac} for all $x$ in $\Ouvs$ and all non-zero $v$. Now take $x$ in $\Ouvs$ and a non-zero vector $v$ and suppose that there exists $T>0$ such that for all $t$ in $[0,T)$, $X_u(x,t)$ is in $\Ouvs$ and $w(t)=\frac{\partial h}{\partial x}(X_u(x,t)) V_u((x,v),t)=0$. To simplify the notations, we denote $X(t)=X_u(x,t)$ and $V(t)=V_u((x,v),t)$. For all integer $i$, we denote 
	\\[0.7em]\null\hfill$\displaystyle 
	w_i(t)=\frac{\partial L_f^{i-1}h}{\partial x}(X(t)) V(t) \ .
	$\hfill\null\\[0.7em]
	We note that for any function $\overline{h}:\RR^n\to\RR$, we have 
	\\[0.7em]$\displaystyle 
	\dot{\overparen{\frac{\partial \overline{h}}{\partial x}(X(t)) V(t)}} = 
	\frac{\partial L_f\overline{h}}{\partial x}(X(t)) V(t)
	$\hfill \null \\\null \hfill $\displaystyle 
	+\sum_{k=1}^{\du}u_k \frac{\partial L_{g_k}\overline{h}}{\partial x}(X(t)) V(t) \ . 
	$\\[0.7em]
	We deduce for all integer $i$ and all $t$ in $[0,T)$
	\\[0.7em]\null\hfill$\displaystyle 
	\dot w_i(t) = w_{i+1}(t)+\sum_{k=1}^{\du}u_k \frac{\partial L_{g_k}L_f^{i-1}h}{\partial x}(X(t)) V(t) \ .
	$\hfill\null\\[0.7em]
	Let us show by induction that $w_i(t)=0$ for all integer $i$ and all $t$ in $[0,T)$. It is true  for $i=1$ by assumption. Now, take an integer $i>1$, and suppose $w_j(t)=0$ for all $t$ in $[0,T)$ and all $j\leq i$, i-e $\frac{\partial \H_i}{\partial x}(X_u(x,t)) V_u((x,v),t)=0$ for all $t$. In particular, $\dot w_i(t)=0$ for all $t$. Besides, according to \eqref{eq_condLipsJac}, $\frac{\partial L_{g_k}L_f^{i-1}h}{\partial x}(X_u(x,t))V_u((x,v),t)=0$ for all $k$ in $\{1,\ldots,\du\}$ and for all $t$. Thus, $w_{i+1}(t)=0$ for all $t$.  
	We conclude that $w_i$ is zero on $[0,T[$ for all $i$ and in particular at time $0$, $\frac{\partial \H_{\od}}{\partial x}(x) v = (w_1(0),\ldots, w_\od(0))=0$. But $\H_{\od}$ is an immersion on $\Ouvs$, thus, necessarily $v=0$ and we have a contradiction. \hfill $\blacksquare$
\end{pf}

\begin{example}
	We go on with Example \ref{ex2}. The linearization of the dynamics \eqref{eq_dyn_ex2} yields
	\\[0.5em]\null \hfill $\displaystyle 
	\dot{v}_1 \;=\;  v_2 \ , \ 
	\dot{v}_2 \;=\; x_3^3v_1 + 3x_3^2x_1v_3 \ , \ 
	\dot{v}_3 \;=\; 0 \ , \
	w \;=\; v_1 \
	\\\null\hfill \refstepcounter{equation}(\arabic{equation}) \label{eq_lin_dyn_ex2}\\$
	Consider $x_0=(x_1,x_2,0)$ in $\Ouvs$ and $v_0=(0,0,v_3)$ with $v_3$ a nonzero real number. 
	The solution to \eqref{eq_dyn_ex2}-\eqref{eq_lin_dyn_ex2} initialized at $(x_0,v_0)$ and with a constant input $u=-1$ is such that $X(x_0,t)$ remains in $\Ouvs$ in $[0,T)$ for some strictly positive $T$ and  $w(t)=0$ for all $t$ in $[0,T)$. Since $v_0$ is nonzero, System \eqref{eq_dyn_ex2} is not uniformly infinitesimally observable on $\Ouvs$. But, for System \eqref{eq_dyn_ex2}, $\H_7$ is an immersion on $\Ouvs$. We deduce from Proposition \ref{prop_CN_infObs} that Inequality \eqref{eq_inequalityLips} is not satisfied for all $i$, i-e there does not exist Lipschitz triangular functions $\gz_i$ for all $i$ on $\Ouvs$.  This is consistent with the conclusion of Example \ref{ex2_suite2}. 
	However, on $\tilde{\Ouvs}$, i-e when we remove the points where $x_3=0$, the system becomes uniformly infinitesimally observable. Indeed, it can be checked that for $x$ in $\tilde{\Ouvs}$, 
	$
	w=\dot{w}=\ddot{w}=w^{(3)}=0
	$, implies necessarily $v=0$.
	Unfortunately, from our results, we cannot infer from this that the functions $\gz_i$ can be taken Lipschitz on $\tilde{\Ouvs}$. Nevertheless, the conclusion of Example  \ref{ex2_suite2} is that $\gz_3$ can be taken Lipschitz even around points with $x_1=0$. All this suggests a possible tighter link between uniform infinitesimal observability and Lipschitzness of the triangular form.
\end{example}

We conclude from this section that uniform infinitesimal observability is required to have the Lipschitzness of the functions $\gz_i$ when they exist. However, we don't know if it is sufficient yet.


\section{Conclusion}

Like for strongly differentially observable systems of order $\dx $ (the system state dimension),
uniform observability of weakly differentially observable systems of order $\od >\dx$, may still imply the existence of  an at least up-to-$\dx+1$-triangular canonical form of dimension $\od$ (see \eqref{3}). 
But
\begin{list}{}{%
		\parskip 0pt plus 0pt minus 0pt%
		\topsep 0.5ex plus 0pt minus 0pt%
		\parsep 0pt plus 0pt minus 0pt%
		\partopsep 0pt plus 0pt minus 0pt%
		\itemsep 0.5ex plus 0pt minus 0pt
		\settowidth{\labelwidth}{-}%
		\setlength{\labelsep}{0.5em}%
		\setlength{\leftmargin}{\labelwidth}%
		\addtolength{\leftmargin}{\labelsep}%
	}
	\item[-]  we have shown this under the additional assumption  that the functions
	$\H_i(x) = (h(x),L_fh(x),...,L_f^{i-1}h(x))$ are open maps. Actually it is sufficient that
	the properties $\B(2),\dots, \B(\dx+1)$ hold (see \eqref{prop_Bi}). 
	\item[-] the functions in the triangular form are possibly non Lipschitz, but only close to points where the rank of the Jacobian of $\H_i$ changes.  Anyhow, uniform infinitesimal observability is necessary to have Lipschitz functions. 
	\item[-] for a non Lipschitz triangular canonical form, convergence of the regular high gain observer may be lost, but, under certain regularity assumptions, it is still possible to design asymptotic observers  (see \cite{BerPraAndAut20162}).
\end{list}

Although our result only gives a partial triangular form, we have no counter example showing that it cannot be a full triangular form.

\supprimer{
	Both uniform infinitesimal observability and uniform observability at least are necessary to have a Lipschitz triangular form which ensures the convergence of a high gain observer. Those conditions being quite restrictive, it is necessary to look for observers which can be used on only continuous triangular forms (see \cite{BerPraAndNOLCOS}).
}

\appendix
\section{Proof of Lemma \ref{claimUniObs}}
\label{appPrLemclaimUniObs}

Assume the system is uniformly observable on $\Ouvs$.
We first show that property $\A(1)$ holds.
Suppose there exists $(x^*_a,x^*_b)$ in $\Ouvs^2$ and $k$ in $\{1,...,\du \}$ such that $x^*_a\neq x^*_b$ and
$$
h(x^*_a)=h(x^*_b) \qquad,\qquad L_{g_k}h(x^*_a)\neq L_{g_k}h(x^*_b) \ .
$$ 
Then, the control law $u$
with all its components zero except its $k$th one which is
$$
u_k\;=\; -\frac{L_fh(x_a)- L_fh(x_b)}{L_{g_k}h(x_a)- L_{g_k}h(x_b)}
\  .
$$
is defined on a neighborhood of $(x_a^*,x_b^*)$. 
The corresponding solutions $X_{u}(x_a^*,t)$ and $X_{u}(x_b^*,t)$ are defined on some time interval $[0,T)$
and satisfy
$$
h(X_{u}(x_a^*,t))=h(X_{u}(x_b^*,t)) \qquad \forall t \in [0,T) \  .
$$
Since $x_a^*$ is different from $x_b^*$,
this contradicts the uniform observability. Thus $\A(1)$ holds.

Let now $i$ in $\{2,\ldots,\dx+1\}$ be such that Property $\B(i)$ holds and $\A(j)$ is satisfied for all $j$ 
in $\{1,\ldots, i-1\}$. 
To establish by contradiction that $\A(i)$ holds, we
assume this is not the case.
This means that there exists $(x_{a,0}^*,x_{b,0}^*)$ in $\Ouvs^2$ and $k$ in
$\{1,...,\du \}$ such that
$\textbf{H}_i(x_{a,0}^*)=\textbf{H}_i(x_{b,0}^*)$
but $L_{g_k}L_f^{i-1}(x_{a,0}^*)\neq L_{g_k}L_f^{i-1}(x_{b,0}^*)$.
This implies $x_{a,0}^*\neq x_{b,0}^*$. By continuity of $L_{g_k}L_f^{i-1}$
and  according to $\B(i)$,
there exists $x_a^*$ (resp $x_b^*$) in $\Ouvs$ sufficiently close to $x_{a,0}^*$ (resp $x_{b,0}^*$) satisfying 
$x_{a}^*\neq x_{b}^*\,  $,
\\[0.5em]\null \hfill $\displaystyle 
\textbf{H}_i(x_{a}^*)=\textbf{H}_i(x_{b}^*)
\  ,\quad 
L_{g_k}L_f^{i-1}(x_{a}^*)\neq L_{g_k}L_f^{i-1}(x_{b}^*)
\: ,
$\hfill \null \\[0.5em]
and $\frac{\partial\textbf{H}_{i-1}}{\partial x}$
is full-rank at $x_a^*$ or $x_b^*$. Without loss of generality, we suppose it is full-rank at $x_a^*$. Thus,
$\frac{\partial\textbf{H}_{j}}{\partial x}$ is full-rank at $x_a^*$ for all
$j<i$. We deduce that there exists an open neighborhood $V_a$ of $x_a^*$ such that for all $j<i$, 
$\frac{\partial\textbf{H}_{j}}{\partial x}$ is full-rank on $V_a$. Since $\A(j)$ holds for all $j<i$, 
according to Lemma \ref{lemRgCst},  $\H_j(V_a)$ is open for all $j<i$ and there exist locally Lipschitz functions $\gz_{j}: \H_j(V_a)\to \RR^{\du}$ such that, for all $x_\alpha$ in 	$V_a$,
\\[0.5em]\null \hfill$\displaystyle
\gz_j(\textbf{H}_j(x_\alpha))=L_gL_f^{j-1}h(x_\alpha) \ .
\ \hfill\null \refstepcounter{equation}(\arabic{equation}) \label{eq_gzVa} $\\[0.5em]
Also, $\textbf{H}_j(x_a^*)=\textbf{H}_j(x_b^*) $ implies that
$\textbf{H}_j(x_b^*) $ is in the open set $\textbf{H}_j(V_a)$. Continuity of each $\textbf{H}_j$ 
implies the existence of an open neighborhood $V_b$ of $x_b^*$ such 
that $\textbf{H}_j(V_b)$ is contained in $\textbf{H}_j(V_a)$ for all
$j<i$.  Thus, for any $x_\beta$ in $V_b$, $\textbf{H}_j(x_\beta)$ is in $\textbf{H}_j(V_a)$, and there exists
$x_\alpha$ in $V_a$ such that $\textbf{H}_j(x_\alpha)=\textbf{H}_j(x_\beta)$. According to $\A(j)$
this implies that $L_gL_f^{j-1}h(x_\beta)=L_gL_f^{j-1}h(x_\alpha)$ and with \eqref{eq_gzVa},
\\[0.4em]$\displaystyle 
L_gL_f^{j-1}h(x_\beta)=L_gL_f^{j-1}h(x_\alpha)=
$\hfill \null \\\null \hfill $\displaystyle 
\gz_j(\textbf{H}_j(x_\alpha))=\gz_j(\textbf{H}_j(x_\beta))\  .
$\\[0.4em]
Therefore, \eqref{eq_gzVa} holds on $V_a$ and $V_b$.

Then, the control law $u$
with all its components zero except its $k$th one which is
\\[0.5em]\null \hfill $\displaystyle 
u_k\;=\; -\frac{L_f^{i}h(x_a)-L_f^{i}h(x_b)}{L_{g_k}L_f^{i-1}h(x_a)-L_{g_k}L_f^{i-1}h(x_b)}
$\hfill\null\\[0.5em]
is defined on a neighborhood of $(x_a^*,x_b^*)$. 
The corresponding solutions $X_{u}(x_a^*,t)$ and $X_{u}(x_b^*,t)$ are defined on some time interval $[0,T)$ where they remain in $V_a$ and $V_b$ respectively.
Let
$Z_a(t)=\textbf{H}_i(X_{u}(x_a^*,t))$, $Z_b(t)=\textbf{H}_i(X_{u}(x_b^*,t))$
and $W(t)=Z_a(t)-Z_b(t)$ on
$[0,T)$.  Since, for all
$j<i$, 
\eqref{eq_gzVa} holds on $V_a$ and $V_b$ , $(W,Z_a)$ is solution
to the system :
\begin{eqnarray*}
	\dot{w}_{1}&\hskip -0.33em=&\hskip -0.33em  w_2 + (\gz_1(z_{a,1})-\gz_1(z_{a,1}-w_1)) \, u \\
	&\hskip -0.33em \ldots&\hskip -0.33em  \\
	\dot{w}_{j}&\hskip -0.33em =&\hskip -0.33em  w_{j+1} \\
	&\hskip -0.33em &\hskip -0.33em
	+ (\gz_j(z_{a,1},...,z_{a,j})-\gz_j(z_{a,1}-w_1,...,z_{a,j}-w_j)) \, u \\
	&\hskip -0.33em \ldots&\hskip -0.33em  \\
	\dot{w}_{i} &\hskip -0.33em =&\hskip -0.33em  0 \\
	\dot{z}_{a,1}&\hskip -0.33em =&\hskip -0.33em  z_2 + \gz_1(z_{a,1}) \, u \\
	&\hskip -0.33em \ldots&\hskip -0.33em  \\
	\dot{z}_{a,j}&\hskip -0.33em =&\hskip -0.33em  z_{j+1} + \gz_j(z_{a,1},...,z_{a,j}) \, u \\
	&\hskip -0.33em \ldots&\hskip -0.33em  \\
	\dot{z}_{a,i}
	&\hskip -0.5em =&\hskip -0.5em  \tilde{u}
\end{eqnarray*}
with initial condition 
$(0,\textbf{H}_i(x_a^*))$,
where $\tilde{u}$ is the time derivative of
$Z_{a,i}(t)$. Note that the function $(0,Z_a)$ is also a solution to this system with the same initial condition. 
Since the functions involved in this system are locally Lipschitz, it admits a unique solution.  Hence, for all $t$ in $[0,T[$, $W(t)=0$, and thus
$Z_a(t)=Z_b(t)$, which implies $h(X(x_a^*,t))=h(X(x_b^*,t))$. Since $x_a^*$ is different from $x_b^*$, this
contradicts the uniform observability. Thus $\A(i )$ holds.


\section{Technical lemmas}
\label{appTechLem}

In this appendix, we consider two continuous functions $\Phi : \RR^n \rightarrow\RR^r$ and $\gamma: \RR^n \rightarrow\RR^q$ and a subset $\Ouvs$ of $\RR^n$ such that 
\begin{equation}
\Phi(x)=\Phi(y) \quad \forall (x,y)\in \Ouvs^2 \, : \, \gamma(x)=\gamma(y) \ .
\label{eqTechLem}
\end{equation}

\begin{lem}
	\label{lem_existenceMap}
	There exists a function $\phi$ defined on $\gamma(\Ouvs)$ such that 
	\begin{equation}
	\label{eq_phi}
	\Phi(x)=\phi(\gamma(x)) \quad \forall x \in \Ouvs \ .
	\end{equation}
\end{lem}

\begin{pf}
	Define the map $\phi$ on $\gamma(\Ouvs)$ as 
	$$
	\phi(z) = \bigcup_{\stackrel{x\in \Ouvs}{\gamma(x)=z}} \{\Phi(x)\} \ .
	$$
	For any $z$ in $\gamma(\Ouvs)$, the set $\phi(z)$ is non-empty and single-valued  because according to \eqref{eqTechLem}, if $z=\gamma(x_a)=\gamma(x_b)$, then $\Phi(x_a)= \Phi(x_b)$. Therefore, we can consider $\phi$ as a function defined on $\gamma(\Ouvs)$ and it verifies \eqref{eq_phi}.
	\hfill $\blacksquare$
\end{pf}

\begin{lem}
	\label{lemTech}
	Consider any compact subset $\C$ of $\Ouvs$. 
	There exists a class $\mathcal{K}$ function $\rho$ such that for all $(x_a,x_b)$ in $\C^2$ 
	\begin{equation}
	|\Phi(x_a)-\Phi(x_b)|\leq \rho\left(|\gamma(x_a)-\gamma(x_b)| \right) \ .
	\label{eqTechLem2}
	\end{equation}
\end{lem}

\begin{pf}
	We denote
	$
	D(x_a,x_b)= |\gamma(x_a)-\gamma(x_b)|
	$.
	Let
	$$
	\rho_0(s) = \max_{\begin{array}{c} (x_a,x_b)\in \C^2\\D(x_a,x_b) \leq s\end{array}}
	|\Phi(x_a)-\Phi(x_b)|
	$$
	This defines properly a non decreasing function 
	with non negative values 
	which satisfies~:
	$$
	|\Phi(x_a)-\Phi(x_b)|\; \leq \; \rho _0(D(x_a,x_b))\qquad \forall (x_a,x_b)\in \C^2
	\  .
	$$
	Also
	$\rho _0(0)=0$. Indeed if not
	there would exist
	$(x_a,x_b)$ in $C^2$ satisfying~:
	$$
	D(x_a,x_b)\;=\; 0\quad ,\qquad |\Phi(x_a)-\Phi(x_b)|\; >\; 0
	\  .
	$$
	But this contradicts Equation (\ref{eqTechLem}).\\
	Moreover, it can be shown that this function is also continuous at $s=0$.
	Indeed, let $(s_k)_{k\in\mathbb N}$ be a sequence converging to $0$. For each $k$, there exist $(x_{a,k}, x_{b,k})$ in $C^2$ which satisfies $D(x_{a,k},x_{b,k})\leq s_k$ and 
	$\rho _0(s_k) = |\Phi(x_{a,k})-\Phi(x_{b,k})|$.
	The sequence $(x_{a,k},x_{b,k})_{k\in\mathbb N}$
	being in a compact set, it
	admits an accumulation point $(x_{a}^*,x_{b}^*)$ which, because of the continuity 
	of $D$ must satisfy 
	$
	D(x_{a}^*, x_{b}^*)=0
	$
	and therefore with (\ref{eqTechLem}) also
	$
	\Phi(x_{a}^*)-\Phi(x_{b}^*)\;=\; 0
	\  .
	$
	It follows that
	$\rho _0(s_k)$ tends to $0$ and $\rho_0$ is continuous at $0$.\\
	Now, the function, defined by the Riemann integral
	$$
	\rho(s) = \left\{
	\begin{array}{lr}
	\frac{1}{s} \displaystyle \int_s^{2s} \rho_0(s) ds + s & , \ s>0 \\
	0 & , \ s=0
	\end{array}
	\right. 
	$$
	is continuous and  strictly increasing and we have~:
	\\[0.7em] \null \hfill$\displaystyle
	|\Phi(x_a)-\Phi(x_b)|\; \leq \; \rho(D(x_a,x_b))
	\qquad 
	\forall (x_a,x_b) \in \C^2\  .
	$
	\hfill $\blacksquare$
\end{pf}

\begin{lem}
	\label{lemTech2}
	Consider any compact subset $\C$ of $\Ouvs$.  There exists
	a continuous function $\phi$ defined on $\RR^q$ such that
	$$
	\Phi(x)=\phi(\gamma(x)) \quad \forall x \in \C \ .
	$$
\end{lem}

\begin{pf}
	Consider $\phi$ and $\rho$ given by Lemmas \ref{lem_existenceMap} and \ref{lemTech} respectively. For any $(z_a,z_b)$ in $\gamma(\C)^2$, there exists $(x_a,x_b)$ in $\C^2$ such that $z_a=\gamma(x_a)$ and $z_b=\gamma(x_b)$. Applying \eqref{eqTechLem2} to $(x_a,x_b)$ and using \eqref{eq_phi}, we have
	$$
	|\phi(z_a)-\phi(z_b)|\; \leq \; \rho(|z_a-z_b|) \ .
	$$
	This means that $\phi$ is uniformly continuous on the compact set $\gamma(\C)$. Thus, $\phi$ is also bounded on $\gamma(\C)$. We deduce from \cite[Corollary 2]{McS} (applied to each component of $\phi$) that $\phi$ admits a uniformly continuous (and bounded) extension defined on $\RR^q$ . 
	\hfill $\blacksquare$
\end{pf}

\begin{lem}
	\label{lemRgCst}
	Assume  that $q\leq n$ and consider an open subset $V$ of $\Ouvs$ such that $\displaystyle \frac{\partial \gamma}{\partial x}$
	is full-rank on $V$, namely $\gamma$ is a submersion on $V$. Then,  $\gamma(V)$ is open and there exists a $C^1$ function $\phi$ defined on $\gamma(V)$  such that 
	$$
	\Phi(x)=\phi(\gamma(x)) \quad \forall x \in V \ .
	$$
\end{lem}

\begin{pf}
	$\gamma$ is an open map according to  \cite[Proposition
	4.28]{Lee}, thus $\gamma(V)$ is open. Consider the function $\phi$ given by Lemma \ref{lem_existenceMap} and take any $z^*$ in $\gamma(V)$. There exists $x^*$ in $V$ such that $z^*=\gamma(x^*)$. $\gamma$ being full-rank at $x^*$, according to the constant rank theorem, there exists an open neighborhood $V^*$ of $x^*$  and $C^1$ diffeomorphisms $\psi_1 : \RR^n\to V^*$ and $\psi_2:\RR^q\to \gamma(V^*)$ such that for all $\tilde{x}$ in $\RR^n$:
	$$
	\gamma(\psi_1(\tilde{x})) = \psi_2(\tilde{x}_1,\ldots,\tilde{x}_q) \ .
	$$
	It follows that for all $z$ in $\gamma(V^*)$ 
	$$
	\gamma(\psi_1(\psi_2^{-1}(z),0))=z
	$$
	namely $\gamma$ admits a $C^1$ right-inverse $\gamma^{ri}$ defined on $\gamma(V^*)$ which is an open neighborhood of $z^*$. 
	Therefore, $\phi=\Phi\circ \gamma^{ri}$ and $\phi$ is $C^1$ at $z^*$. 
	\hfill $\blacksquare$
\end{pf}

\bibliographystyle{elsart-harv}        
\bibliography{biblio}

\end{document}